\documentclass[12pt]{amsart}

\usepackage[english]{babel}
\usepackage{tikz-cd}
\usepackage{amsrefs}
\usepackage{graphicx}
\graphicspath{ {images/} }
\usepackage[T1]{fontenc}
\usepackage{times}
\usepackage{amssymb,verbatim,amscd,amsfonts,amsmath,amsthm,enumerate}
\usepackage[all]{xy}
\usepackage{subfig}
\usepackage{tikz}
\usetikzlibrary{shapes,arrows,shadows}
\usetikzlibrary{decorations.markings}
\usepackage{color}
\usepackage[normalem]{ulem}
\usepackage{hyperref}
\usepackage{wrapfig}
\usetikzlibrary{arrows}
\usepackage{pb-diagram}

\newtheorem{theorem}{Theorem}[section]

\newtheorem{corollary}[theorem]{Corollary}
\newtheorem{lemma}[theorem]{Lemma}

\newtheorem{remark}[theorem]{Remark}

\begin{document}

%%%%%%%%%%%%%%%%%%%%%%%%%%%%%%%%%%%%%%%
%%%%%%%%%%%%%%%%%%%%%%%%%

%%%%%%%%%%%%%%%%%%%%%%%%%%%%%%%%%%%
%%%%%%%%%%%%%%%%%%%%%%%%%%%%%%%%%%%
%%%%%%Authors
%%%%%%%%%%%%%%%%%%%%%%%%%%%%%%%%%%
%%%%%%%%%%%%%%%%%%%%%%%%%%%%%%%%%%
\title{On the Cauchy Integral Theorem and Polish spaces}
\author{Cristian L\'opez Morales and Camilo Ram\'irez Maluendas}

\address{Departamento de Matem\'aticas y Estad\'istica, Universidad Nacional de Colombia, Sede Manizales, Manizales 170004, Colombia.}
\email{crlopezmo@unal.edu.co, camramirezma@unal.edu.co}

\thanks{The second author was partially supported by Proyecto Hermes 59315}
\maketitle

%%%%%%%%%%%%%%%%%%%%%%%%%%%%%%%%%%%%%%%%%%%
%%%%%%%%%%%%%%%%%%%%%%%%%%%%%%%%%%%%%%%%%%%
%%%%%%%%%%%Abstract
%%%%%%%%%%%%%%%%%%%%%%%%%%%%%%%%%%%%%%%%%%%
%%%%%%%%%%%%%%%%%%%%%%%%%%%%%%%%%%%%%%%%%%%

\begin{abstract}
We prove that if a function $f$ is continuous in an open subset $U\subset\mathbb{C}$ and analytic in $U\setminus X$, where $X\subset U$ is a Polish space having characteristic system $(i,n)$, such that $i\in\{0,1\}$ and $n\in\mathbb{N}$, then the complex integral line of $f$ along the boundary of any triangle in $U$ vanishes. 
\end{abstract}

%%%%%%%%%%%%%%%%%%%%%%%%%%%%%%%%%%%%%%
%%%%%%%%%%%%%%%%%%%%%%%%%%%%%%%%%%%%%%
%%%%%%%%%%%Introduction
%%%%%%%%%%%%%%%%%%%%%%%%%%%%%%%%%%%%%%
%%%%%%%%%%%%%%%%%%%%%%%%%%%%%%%%%%%%%%

\section{Introduction}

In 1814, Augustin Louis Cauchy presented to the \emph{French Academy of Sciences} the document \emph{M\'emoire sur les integrales d\'efinies} containing his contributions to development of the theory of complex functions \cite{Ett}. In his memoirs, Cauchy proved that if a complex function $f$ is analytic within a closed curve $\gamma$  and also on the curve itself, then the complex integral line of $f$ around that curve, is equal to zero. In the two works \cites{Goursat1884, Goursat1900}, Edouard Goursat proved Cauchy Integral Theorem:
\[
\int\limits_{\gamma}f(z)dz=\boldsymbol{0},
\]
without assuming the continuity of the derivative $f'(z)$ on the closed region $U$ bounded by the curve of integration $\gamma$. For a short summary of these works, refer to \cite{Col1902}*{p.427--429}. In 1900, Eliakim Hastings Moore   wrote and published his proof of the Cauchy Integral Theorem \cite{More1900}. One year later, Alfred Pringsheim  also wrote his version of the proof \cite{Pring1901}. For more details about of the historical development of its proof see \cite{Scott1978}. 

Motivated by the Cauchy Integral Theorem and the several research of its proof, one might ask: \emph{What are the weakest set of the assumptions for  what the complex integral line of an "analytic function" along a closed curve vanishes?}

Recall that a \emph{Polish space} is a separable completely metrizable topological space. Given a Polish space $X$, the \emph{Cantor-Bendixson derivative of $X$} is the set 
\[
X':= \{x \in X : x \text{ is a limit point of } X\}.
\]

For every $k\in\mathbb{N}$, the \emph{$k$-th Cantor-Bendixson derivative of $X$} is defined as  $X^{k}:=(X^{k-1})'$, where $X^{0} = X$ (for details see \cite{kechris}*{p. 34}). A countable Polish space is said to have \textit{characteristic system} $(k, n)$ if $X^{k}\neq \emptyset$, $X^{k+1} = \emptyset$ and the cardinality of $X^k$ is $n$, for some natural number $n\in\mathbb{N}$.

Any countable Polish space with characteristic system $(k, n)$ is homeomorphic to the ordinal number $\omega^{k}
+n$, where $\omega$ is the least infinite ordinal. In particular, every two countable Polish spaces with characteristic system $(k, n)$ are homeomorphic. For further details, see \cite{MaSi}. We remark that any finite Polish space (finite discrete space) of cardinality $n$, has characteristic system $(0,n)$, for some $n\in\mathbb{N}$. In addition, the closed subset $X=\overline{\left\{\frac{1}{n}:n\in\mathbb{N}\right\}}\subset\mathbb{R}$ has characteristic system $(1,1)$. In other words, this compact space $X$ is homeomorphic to the ordinal number $\omega+1$.

In the classical case, when $U$ is an open subset of $\mathbb{C}$, $\triangle(z_{1},z_{2},z_{3})\subset U$ is a triangle with vertices $z_{1}$, $z_{2}$ and $z_{3}$, and $f$ is an analytic function in $U$; then the Cauchy Integral Theorem implies
\[
\int\limits_{[z_{1},z_{2},z_{3},z_{1}]}f(z)dz=\boldsymbol{0},
\]
when $[z_{1},z_{2},z_{3},z_{1}]$ is the polygonal closed curve parametrizing the boundary $\partial \triangle(z_{1},z_{2},z_{3})$ of the triangle $\triangle(z_{1},z_{2},z_{3})$. It is possible to replace the assumption $f$ analytic in whole $U$ by a weak one. More precisely,

\begin{lemma}[\cite{Palka}*{Lemma 1.2. p. 143}]\label{twopoints}
If a function $f$ is continuous in an open subset $U\subset\mathbb{C}$ and analytic in $U\setminus\{w\}$ for some point $w$ of $U$, then

\[
\int\limits_{[z_1,z_2,z_3,z_1]}f(z) dz=\boldsymbol{0},
\]
for every  triangle $\triangle(z_1,z_2,z_3)\subset U$.
\end{lemma}

We extend the above Lemma for the case when $f$ is analytic in $U$ removing finitely many points. More precisely, we prove the following result:

\begin{theorem}\label{t:n_points}
If a function $f$ is continuous in an open subset $U\subset\mathbb{C}$ and analytic in $U\setminus X$, where $X\subset U$ is a finite Polish space having characteristic system $(0,n)$, for some $n\in\mathbb{N}$. Then
\[
\int\limits_{[z_1,z_2,z_3,z_1]}f(z)dz=\boldsymbol{0},
\]
for every  triangle $\triangle(z_1,z_2,z_3)\subset U$.
\end{theorem}

Even more, we provide a more general form to the previous Theorem, it means, the complex integral line also is zero, when the Polish space $X$ has characteristic system $(i,n)$, when $i\in\{0,1\}$ and $n\in\mathbb{N}$. Our mains contributions are as follows:

\begin{theorem}\label{theorem:characteristic_system_(1,1)}
If a function $f$ is continuous in an open subset $U\subset \mathbb{C}$ and analytic in $U\setminus X$, where $X\subset U$ is a Polish space having characteristic system $(1,n)$, for some $n\in\mathbb{N}$. Then
\[
\int\limits_{[z_{1},z_{2},z_{3},z_{1}]} f(z)dz=\boldsymbol{0},
\]
for every triangle $\triangle(z_{1},z_{2},z_{3})\subset U$.
\end{theorem}

From the previous Theorems, it is immediate the following corollary:

\begin{corollary}
If a function $f$ is continuous in an open subset $U\subset\mathbb{C}$ and analytic in $U\setminus X$, where $X\subset U$ is a Polish space having characteristic system $(i,n)$, such that $i\in\{0,1\}$ and $n\in\mathbb{N}$, then the complex integral line of $f$ along the boundary of any triangle in $U$ vanishes. 
\end{corollary}

\subsection*{Acknowledgments.} The authors were partially supported by \textbf{Universidad Nacional de Colombia, Sede Manizales}. The first author has dedicated this work in memory of his grandfather Jaime, who always believed in him. The second author has dedicated this work to his beautiful family: Marbella and Emilio, in appreciation of their love and support.

%%%%%%%%%%%%%%%%%%%%%%%%%%%%%%%%%%%%%%%%
%%%%%%%%%%%%%%%%%%%%%%%%%%5
\section{Proof of the Theorem \ref{t:n_points}}

We consider $f$ a continuous map in the open subset $U\subset \mathbb{C}$ and analytic in $U\setminus X$, such that $X\subset U$ is a finite Polish space. We shall proceed by induction on the cardinality of $X$. 

\subsection*{The Polish space $X$ has cardinality one.} In other words, $X$ has characteristic system $(0,1)$. From the Lemma \ref{twopoints} we hold that the complex integral  line of $f$ along the boundary of any triangle in $U$ is zero. 

\subsection*{The Polish space $X$ has cardinality $k$, for each $k<n$ and some natural number fixed $n\in\mathbb{N}$.} It implies that $X$ has characteristic system $(0,k)$. Then, it satisfies that
\begin{equation*}\label{eq:the1_induction_1}
\int\limits_{[z_1,z_2,z_3,z_1]} f(z)dz=\boldsymbol{0},
\end{equation*}
for any triangle $\triangle(z_1,z_2,z_3)\subset U$. 

\subsection*{The Polish space $X$ has cardinality $n$.} 
\begin{figure}[ht]
	\begin{center}
		\begin{tikzpicture}[baseline=(current bounding box.north)]
		\begin{scope}[scale=1.0]
		\clip (-3.7,-0.2) rectangle (4.7,3.7);
		%%%%%%%%%%%%%%%%%%%%%%%%%%%%%%%%%%%%%%%%%%%%%%%
		%%%%%   
		%%%%%%%%%%%%%%%%%%%%%%%%%
		%%%%%%%%%%%%%%%%%%%%%%%%%%
             \draw [draw=blue!50, fill=blue!10, line width=1pt] (-2,0.3) -- (2.2,1.7) -- (-0.8,3.4) -- cycle;
		%%%%%%%%%Plano complejo
  \draw [line width=1.0pt, black!50] (-0.8,0.7) -- (-0.2,2.8);
		\draw [line width=1pt] (-3.5,0) -- (3.5,0);
		\draw [line width=1pt] (-3.5,0) -- (-2.5,3.5);  
		\draw [line width=1pt] (3.5,0) -- (4.5,3.5); 	
		\draw [line width=1pt] (-2.5,3.5) -- (4.5,3.5);	
		%%%%%%%%%%%%%%%%%%%%%
		\draw [line width=1pt, red] (-3,1.7) -- (4,1.7);
		%%%%%%%%%%%%etiquetas
		\node at (-0.7,2.8) {$w_{j}$};
		\node at (-0.2,2.8) {$\bullet$};
             \node at (-1.2,0.7) {$w_{i}$};
		\node at (-0.8,0.7) {$\bullet$};
		\node at (2.2,1.7) {$\bullet$};
		\node at (2.2,2.1) {$z_{1}$};
            \node at (-0.53,1.7) {$\bullet$};
            \node at (-0.7,2.1) {$v$};
            \node at (-1.45,1.7) {$\bullet$};
            \node at (-1.6,2.1) {$w$};
		\node at (-0.4,1.2) {$\ell$};
            \node at (3,2.1) {$\ell'$};
		\node at (-2.8,0.5) {$\mathbb{C}$};
		\node at (2,0.5) {$C_{1}$};
            \node at (2.5,3) {$C_{2}$};
		\end{scope}
		\end{tikzpicture}
	\end{center}
	\caption{\emph{Straight lines $\ell$ and $\ell'$.}}
	\label{Fig:straight_line_1}
\end{figure}
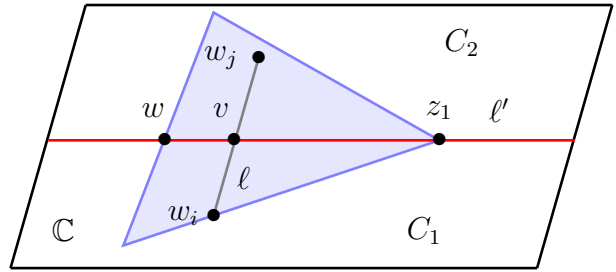

It means, $X$ has characteristic system $(0,n)$. We denote by $w_{1},\ldots,w_{n}$ the $n$ points in $U$ such that $X=\{w_{1},\dots,w_{n}\}$. We shall prove that
\[
\int\limits_{[z_{1},z_{2},z_{3},z_{1}]}f(z)dz=\boldsymbol{0},
\]
for any triangle $\triangle(z_{1},z_{2},z_{3})\subset U$. The induction hypothesis implies that the above complex integral line is equal to zero, when at least one of the these points $w_{1},\ldots, w_{n}$ lies outside of $\triangle(z_{1},z_{2},z_{3})$. Now, we shall study the case when all points $w_{1},\ldots,w_{n}$ are in $\triangle(z_{1},z_{2},z_{3})$. Then, we choose two of these points $w_{i}$ and $w_{j}$ and draw the straight line $\ell$ passing through them. Since the three points $z_{1},z_{2}$ and $z_{3}$ are not collinear, then we can suppose without loss of generality that $z_{1}$ is not in $\ell$. Thus, we draw the straight line $\ell'$ passing trough $z_{1}$ and the middle point $v$ of the straight line segment with endpoints $w_{i}$ and $w_{j}$, see the Figure \ref{Fig:straight_line_1}.

\begin{remark}\label{remark:Properties_sl_polish.space}
The straight line $\ell'$ satisfies the following properties:

\begin{itemize}
\item[\textbf{(a)}] It decomposes the complex plane $\mathbb{C}$ into two open connected subsets $C_{1}$ and $C_{2}$, such that the vertex $z_{2}$ is belonged to one of these open connected, and the vertex $z_{3}$ is belonged to the other open connected.

\item[\textbf{(b)}] If we choose the point $v_{i}\in C_{i}$, with $i\in\{1,2\}$, then the straight line segment having endpoints $v_{1}$ and $v_{2}$, intersects the straight line $\ell'$ in a unique point. In particular, the points $w_{i}$ and $w_{j}$ are not belonged to $\ell'$.
\end{itemize}
\end{remark}

The paragraphs \textbf{(a)} and \textbf{(b)} above described imply that $\ell'$ must intersect the straight line segment with endpoints $z_{2}$ and $z_{3}$, the unique common point is denoted by $w$. Thus, we decompose $\triangle(z_{1},z_{2},z_{3})$ into the triangles $\triangle(z_{1},z_{2},w)$ and $\triangle(z_{1},w,z_{3})$, which satisfy the following properties:
\begin{itemize}
    \item[\textbf{(1)}] They have in common the straight line segment $L$ with endpoints $z_{1}$ and $w$. Moreover, $L\subset \ell'$.
    \item[\textbf{(2)}] Each one of these triangles $\triangle(z_1,z_2,w)$ and $\triangle(z_1,w,z_3)$, have less than $n$ points of $X$.
\end{itemize}
Therefore, we have the following equality between complex integral line
\begin{equation}\label{Eq:XXX2}
\int\limits_{[z_{1},z_{2},z_{3},z_{1}]} f(z) dz=\int\limits_{[z_{1},z_{2},w,z_{1}]} f(z) dz+\int\limits_{[z_{1},w,z_{3},z_{1}]}f(z)dz.
\end{equation}
By the property \textbf{(2)} and the induction hypothesis, we hold that for $\triangle(z_{1},z_{2},w)$ and $\triangle(z_{1},w,z_{3})$ it satisfies that 
\begin{equation}\label{Eq:Induction}
    \int\limits_{[z_{1},z_{2},w,z_{1}]} f(z) dz=\int\limits_{[z_{1},w,z_{3},z_{1}]}f(z)dz=\boldsymbol{0}.
\end{equation}
Now, replacing the values of equation \eqref{Eq:Induction} into the equation \eqref{Eq:XXX2}, we conclude that
\[
\int\limits_{[z_{1},z_{2},z_{3},z_{1}]} f(z) dz=\boldsymbol{0}.
\]
\qed

%%%%%%%%%%%%%%%%%%%%%%%%%%%%%%%%%%%%%%%%%%%%%%%%%%%%%%%%%%%%%%%%%%%%%%%%%%%%%
\section{Proof of the Theorem \ref{theorem:characteristic_system_(1,1)}}

We shall develop same ideas appearing in the proof of the Theorem \ref{t:n_points}. Let $f$ be a continuous map in the open subset $U\subset \mathbb{C}$ and analytic in $U\setminus X$, such that $X\subset U$ is a Polish space. We take a triangle $\triangle(z_{1},z_{2},z_{3})\subset U$. We shall proceed by induction on the cardinality of  the $1$-st Cantor-Bendixson derivative $X'$ of $X$.

\subsection*{The set $X'$ has only one point.} It means, the characteristic system of $X$ is $(1,1)$. We denote by $w$ the only one point of $X'$. We shall study the following cases.

\subsubsection*{Case 1. The point $w$ lies outside of $\triangle(z_{1},z_{2},z_{3})$} Then we can found an open subset $V\subset U$, such that $\triangle(z_{1},z_{2},z_{3}) \subset V$ and $w\notin V$. It implies that the open $V$ has at most a finitely many points of $X$. Applying the Theorem \ref{t:n_points} on the open $V$, we obtain the expected value  for the complex integral line. More precisely,
\[
\int\limits_{[z_{1},z_{2},z_{3},z_{1}]}f(z)dz=\boldsymbol{0}.
\]

\subsubsection*{Case 2. The point $w$ is a vertex of $\triangle(z_{1},z_{2},z_{3})$} 
\begin{figure}[ht]
	\centering
	\begin{tabular}{ccccc}	
	\begin{tikzpicture}[baseline=(current bounding box.north)] 
	\begin{scope}[scale=0.7]
	\clip (-5,-0.9) rectangle (3,4.5);
%%%%%%%%%%%%%%%%%%%%%%%%%%%%%%%%%%%%%%%%%%%%%%%
	\draw [draw=blue!50, fill=blue!10, line width=1.0pt] (-3,0) -- (2,0) -- (1,4) -- cycle;
    \draw [red, line width=1.0pt] (0,0) -- (-1,2);
    \draw [red, line width=1.0pt] (0,0) -- (1,4);
	\node at (-4,0){$w=z_{1}$};
    \node at (0,-0.5){$\boldsymbol{w}_{2}(t)$};
    \node at (0,0){$\bullet$};
    \node at (-1.9,2.2){$\boldsymbol{w}_{3}(t)$};
    \node at (-1,2){$\bullet$};
	\node at (1,4.3){$z_{3}$};
	\node at (2.4,0){$z_{2}$};
	\end{scope}
	\end{tikzpicture}
  &		
		\begin{tikzpicture}[baseline=(current bounding box.north)] 
	\begin{scope}[scale=0.7]
	\clip (-4.8,-0.9) rectangle (2.6,4.5);
%%%%%%%%%%%%%%%%%%%%%%%%%%%%%%%%%%%%%%%%%%%%%%%
	\draw [draw=blue!50, fill=blue!10, line width=1.0pt] (-3,0) -- (2,0) -- (1,4) -- cycle;
    \draw [red, line width=1.0pt] (0,0) -- (1,4);
    \node at (0,-0.5){$w$};
    \node at (0,0){$\bullet$};
	\node at (-3.3,0){$z_{1}$};
	\node at (1,4.3){$z_{3}$};
	\node at (2.4,0){$z_{2}$};
	\end{scope}
	\end{tikzpicture}	
 &
 \begin{tikzpicture}[baseline=(current bounding box.north)] 
	\begin{scope}[scale=0.7]
	\clip (-4.8,-0.9) rectangle (2.6,4.5);
%%%%%%%%%%%%%%%%%%%%%%%%%%%%%%%%%%%%%%%%%%%%%%%
	\draw [draw=blue!50, fill=blue!10, line width=1.0pt] (-3,0) -- (2,0) -- (1,4) -- cycle;
    \draw [red, line width=1.0pt] (-3,0) -- (1.6,1.7);
    \node at (0,1.5){$w$};
    \node at (0,1.1){$\bullet$};
	\node at (-3.3,0){$z_{1}$};
	\node at (1,4.3){$z_{3}$};
	\node at (2.4,0){$z_{2}$};
	\end{scope}
	\end{tikzpicture}	
	\\
		\textbf{(a)} \emph{The point $w$ is a vertex.}                         &  \textbf{(b)} \emph{The point $w$ is in an edge.} 
                   & \textbf{(c)} \emph{The point $w$ is in the interior.}\\
	\end{tabular}
	\caption{\emph{Triangle $\triangle(z_{1},z_{2},z_{3})$.}}
	\label{Fig:def_triangle_20}
\end{figure}
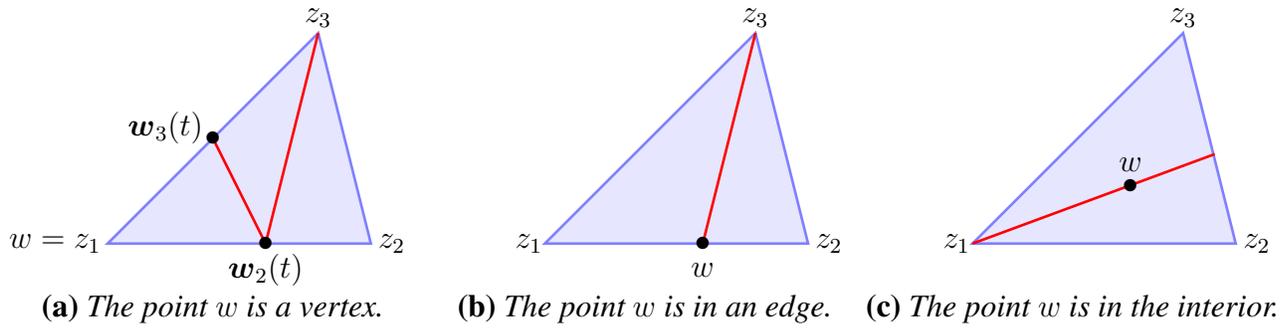

We can suppose without loss of generality that $w=z_{1}$. Now, we consider the continuous maps $\boldsymbol{w}_{2}$, $\boldsymbol{w}_{3}:(0,1)\to \triangle(z_{1},z_{2},z_{3})$ given by $\boldsymbol{w}_{2}(t)=(1-t)z_{1}+tz_{2}$ and $\boldsymbol{w}_{3}(t)=(1-t)z_{1}+tz_{3}$, respectively. By construction the three points $\boldsymbol{w}_{2}(t)$, $\boldsymbol{w}_{3}(t)$ and $z_{3}$ (and the three points $\boldsymbol{w}_{2}(t)$, $z_{2}$ and $z_{3}$, respectively) are not collinear as it is shown in the Figure \ref{Fig:def_triangle_20}-\textbf{(a)}. If we consider the decomposition of $\triangle(z_{1},z_{2},z_{3})$ into the triangles $\triangle(z_{1},\boldsymbol{w}_{2}(t),\boldsymbol{w}_{3}(t))$, $\triangle(\boldsymbol{w}_{2}(t),\boldsymbol{w}_{3}(t),z_{3})$, $\triangle(\boldsymbol{w}_{2}(t),z_{2},z_{3})$ and apply the properties of the complex integral line, then we obtain 

\begin{equation}\label{eq:trianles_X_one_point}
\int\limits_{[z_{1},z_{2},z_{3},z_{1}]} f(z)dz=\int\limits_{[z_{1},\boldsymbol{w}_{2}(t),\boldsymbol{w}_{3}(t),z_{1}]}f(z)dz+\int\limits_{[\boldsymbol{w}_{3}(t),\boldsymbol{w}_{2}(t),z_{3},\boldsymbol{w}_{3}(t)]}f(z)dz+\int\limits_{[\boldsymbol{w}_{2}(t),z_{2},z_{3},\boldsymbol{w}_{2}(t)]}f(z)dz.
\end{equation}
We note that $\triangle(\boldsymbol{w}_{2}(t),\boldsymbol{w}_{3}(t),z_{3})$ and $\triangle(\boldsymbol{w}_{2}(t),z_{2},z_{3})$ satisfy the conditions describes in the previous \emph{Case 1}, then the two complex integral line of the right side of this last equality are equal to zero. Then the equation \eqref{eq:trianles_X_one_point} is re-written as follows
\begin{equation}\label{eq:trianles_X_one_point_2}
\int\limits_{[z_{1},z_{2},z_{3},z_{1}]} f(z)dz=\int\limits_{[z_{1},\boldsymbol{w}_{2}(t),\boldsymbol{w}_{3}(t),z_{1}]}f(z)dz.
\end{equation}
Now, we take the modulus of the complex integral line in the previous equation (\ref{eq:trianles_X_one_point_2}) and consider the upper bound
\begin{equation}\label{eq:trianles_X_one_point_3}
\left\vert \, \, \int\limits_{[z_{1},z_{2},z_{3},z_{1}]} f(z)dz\right\vert =\left\vert \,\,  \int\limits_{[z_{1},\boldsymbol{w}_{2}(t),\boldsymbol{w}_{3}(t),z_{1}]}f(z)dz\right\vert \leq m \, \ell,
\end{equation}
where $m=\max\{|f(w)|: w\in \triangle(z_{1},z_{2},z_{3})\}$ and $\ell$ is the length of the boundary of $\triangle(z_{1},\boldsymbol{w}_{2}(t),\boldsymbol{w}_{3}(t))$. It is easy to verify that
$\ell=t(|z_{1}-z_{2}|+|z_{2}-z_{3}|+|z_{1}-z_{3}|)$. Substituting the value of $\ell$ in the equation (\ref{eq:trianles_X_one_point_3}) we obtain
\begin{equation}
\left\vert \, \, \int\limits_{[z_{1},z_{2},z_{3},z_{1}]} f(z)dz\right|\leq m \, t(|z_{1}-z_{2}|+|z_{2}-z_{3}|+|z_{1}-z_{3}|).
\end{equation}
As the above equation is defined for each $t\in(0,1)$, then the right-hand side of this inequality tends to $\boldsymbol{0}$ as $t$ goes to $0$. Thus, we conclude
\[
\int\limits_{[z_{1},z_{2},z_{3},z_{1}]} f(z)dz=\boldsymbol{0}.
\]

\subsubsection*{Case 3. The point $w$ is an edge of $\triangle(z_{1},z_{2},z_{3})$} We can suppose without loss of generality that $w$ is in the edge having endpoints $z_{1}$ and $z_{2}$. Thus, we consider the decomposition of $\triangle(z_{1},z_{2},z_{3})$ into the triangles $\triangle(w,z_{2},z_{3})$ and $\triangle(w,z_{3},z_{1})$ as shown the Figure \ref{Fig:def_triangle_20}-\textbf{(b)}. Then we hold
\begin{equation*}    \int\limits_{[z_{1},z_{2},z_{3},z_{1}]}f(z)dz=\int\limits_{[z_{1},w,z_{3},z_{1}]}f(z)dz+\int\limits_{[w,z_{2},z_{3},w]}f(z)dz.
\end{equation*}
We remark that the both complex integral line of the right side of this last equality are equal to zero, because they are of the form described in the previous \textit{Case 2}. Thus, we conclude
\[
\int\limits_{[z_{1},z_{2},z_{3},z_{1}]}f(z)dz=\boldsymbol{0}.
\]

\subsubsection*{Case 4. The point $w$ is in the interior set of $\triangle(z_{1},z_{2},z_{3})$} Then we consider the straight line segment such that one of its endpoints is $z_{1}$ and the other endpoints is in the opposite edge of this vertex. Furthermore, such straight line segment passes through the point $w$ as shown the Figure \ref{Fig:def_triangle_20}-\textbf{(c)}. Then, we decompose $\triangle(z_{1},z_{2},z_{3})$ into two triangles, which are in the previous \textit{Case 3}. Thus, we get that the value of the complex integral line in the Theorem is zero.

\subsection*{The set $X'$ has cardinality $k$, for $k<n$ and some natural number fixed $n\in\mathbb{N}$.} In other words, the characteristic system of $X$ is $(1,k)$. It implies $X$ has characteristic system $(0,k)$. Then, it satisfies
\begin{equation}
\int\limits_{[z_1,z_2,z_3,z_1]} f(z)dz=\boldsymbol{0}.
\end{equation}

\subsection*{The set $X'$ has cardinality $n$.} 
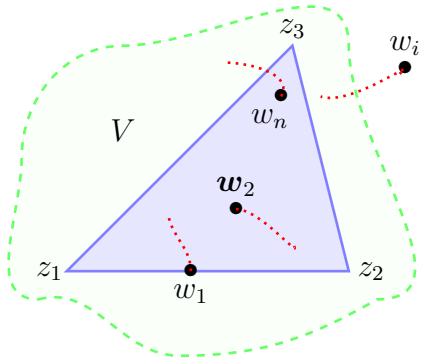
\begin{figure}[ht]
        \centering
        \begin{tikzpicture}[baseline=(current bounding box.north)] 
        \begin{scope}[scale=0.75]
	\clip (-4.8,-2) rectangle (4,5);
    %%%%%%%%%%%%%%%%%%%%%%%%%%%%%%%%%%%%%%%%%%%%%%%
        
        \draw[green!60, fill=green!2, dashed, line width=1pt] (-3.5,-0.5) to[out=15, in=180] (-1,-1.5) to[out=0, in=210] (2.5,-0.5) to[out=0, in=-90] (2,4) to[out=90, in=0] (-2,4) to[out=180, in=90] (-3.8,2) to[out=-90, in=180] cycle;
    
        \draw [draw=blue!50, fill=blue!10, line width=1.0pt] (-3,0) -- (2,0) -- (1,4) -- cycle;
        \node at (-0.8,-0.4){$w_{1}$};
        \node at (-0.8,0){$\bullet$};
        \draw[dotted, red, line width=1.0pt, decorate] (-0.8,0) to[out=90, in=-80] (-1.2,1);
        \node at (0,1.5){$\boldsymbol{w}_{2}$};
        \node at (0,1.1){$\bullet$};
        \draw[dotted, red, line width=1.0pt, decorate] (0,1.1) to[out=0, in=-60] (1,0.5);
        \node at (0.6,2.7){$w_{n}$};
        \node at (0.8,3.1){$\bullet$};
        \draw[dotted, red, line width=1.0pt, decorate] (0.8,3.1) to[out=60, in=0] (-0.2,3.7);
	\node at (-3.3,0){$z_{1}$};
	\node at (1,4.3){$z_{3}$};
	\node at (2.4,0){$z_{2}$};
        \node at (-2,2.5){$V$};
        \node at (3,4){$w_i$};
        \node at (3,3.6){$\bullet$};
        \draw[dotted, red, line width=1.0pt, decorate] (2.9,3.5) to[out=40, in=-20] (1.5,3.1);
    	\end{scope}
    \end{tikzpicture}
    \caption{\emph{Open set $V$ contains the triangle $\triangle(z_{1},z_{2},z_{3})$.}}
    \label{fig:limitpoint}
\end{figure}

In other words, the Polish space $X$ has characteristic system $(1,n)$. We assume that $X'=\{w_1,\dots,w_n\}$. If $w_i$ lies outside of $\triangle(z_1,z_2,z_3)$ for some $i\in\{1,\dots,n\}$, then we define the finite set
    \[
    \widetilde{X}:=X^{'}\cap (U\setminus \triangle(z_{1},z_{2},z_{3})),
    \]
    which has cardinality $k$ such that $0<k<n$.    Thus, there exists an open subset $V\subset U$ such that $\triangle(z_{1},z_{2},z_{3})\subset V$ and $\widetilde{X}\cap V=\emptyset$ (see Figure \ref{fig:limitpoint}). We remark that $Y:=X\cap \triangle(z_{1},z_{2},z_{3})$ is also a Polish space, which has characteristic system $(1,n-k)$. 
    Given that the function $f$ is continuous in $V\subset U$ and analytic in $V\setminus Y$, by the induction hypothesis we get that the complex integral line of $f$ along the polygonal closed curve $[z_{1},z_{2},z_{3},z_{1}]$, is equal to zero.

\begin{figure}[ht]
	\begin{center}
		\begin{tikzpicture}[baseline=(current bounding box.north)]
		\begin{scope}[scale=1.0]
		\clip (-3.7,-0.2) rectangle (4.7,3.7);
		%%%%%%%%%%%%%%%%%%%%%%%%%%%%%%%%%%%%%%%%%%%%%%%
		%%%%%   
		%%%%%%%%%%%%%%%%%%%%%%%%%
		%%%%%%%%%%%%%%%%%%%%%%%%%
             \draw [draw=blue!50, fill=blue!10, line width=1pt] (-2,0.3) -- (2.2,1.7) -- (-0.8,3.4) -- cycle;
		%%%%%%%%%Plano complejo
  \draw [line width=1.0pt, black!50] (-0.8,0.7) -- (-0.2,2.8);
		\draw [line width=1pt] (-3.5,0) -- (3.5,0);
		\draw [line width=1pt] (-3.5,0) -- (-2.5,3.5);  
		\draw [line width=1pt] (3.5,0) -- (4.5,3.5); 	
		\draw [line width=1pt] (-2.5,3.5) -- (4.5,3.5);	
		%%%%%%%%%%%%%%%%%%%%%
		\draw [line width=1pt, red] (-3,1.7) -- (4,1.7);
		%%%%%%%%%%%%etiquetas
		\node at (-0.6,2.8) {$w_{j}$};
		\node at (-0.2,2.8) {$\bullet$};
            \draw[dotted, red, line width=1.0pt, decorate] (-0.2,2.8) to[out=-20, in=0] (0.3,2.9);
            \node at (-1.2,0.9) {$w_{i}$};
		\node at (-0.8,0.7) {$\bullet$};
            \draw[dotted, red, line width=1.0pt, decorate] (-0.8,0.7) to[out=0, in=-60] (-0.3,0.5);
		\node at (2.2,1.7) {$\bullet$};
		\node at (2.2,2.1) {$z_{1}$};
            \node at (-0.55,1.7) {$\bullet$};
            \node at (-0.7,2.1) {$v$};
            \node at (-1.45,1.7) {$\bullet$};
            \node at (-1.6,2.1) {$w$};
		\node at (-0.4,1.2) {$\ell$};
            \node at (3,2.1) {$\ell'$};
		\node at (-2.8,0.5) {$\mathbb{C}$};
		\node at (2,0.5) {$C_{1}$};
            \node at (2.5,3) {$C_{2}$};
		\end{scope}
		\end{tikzpicture}
	\end{center}
	\caption{\emph{Connected component $C_{1}$ and $C_{2}$}.}
	\label{Fig:straight_line_3}
\end{figure}
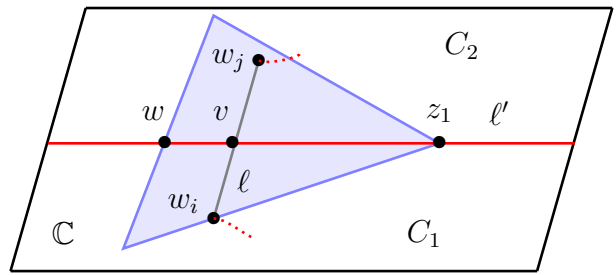

Now, we suppose that $X^{'}\subset\triangle(z_1,z_2,z_3)$. Then, we take the straight line $\ell$ passing through the different points $w_{i}$ and $w_{j}$ in $X'$. Since the three points $z_{1},z_{2}$ and $z_{3}$ are not collinear, then we can suppose without loss of generality that the point $z_{1}$ is not in the straight line $\ell$. Thus, we draw the straight line $\ell'$ passing trough $z_{1}$ and the middle point $v$ of the straight line segment with endpoints $w_{i}$ and $w_{j}$ as shown the Figure \ref{Fig:straight_line_3}.

We note that the straight line $\ell'$ satisfies the same properties described in the Remark \ref{remark:Properties_sl_polish.space}. The paragraphs \textbf{(a)} and \textbf{(b)} described in such Remark, imply that $\ell'$ must intersect the straight line segment with endpoints $z_{2}$ and $z_{3}$, the unique common point is denoted by $w$. Thus, we decompose $\triangle(z_{1},z_{2},z_{3})$ into the triangles $\triangle(z_{1},z_{2},w)$ and $\triangle(z_{1},w,z_{3})$. Then we have the following equality
\begin{equation}\label{eq:final}
\int\limits_{[z_{1},z_{2},z_{3},z_{1}]} f(z) dz=\int\limits_{[z_{1},z_{2},w,z_{1}]} f(z) dz+\int\limits_{[z_{1},w,z_{3},z_{1}]}f(z)dz.
\end{equation}

On the other hand, we remark that $Y_{1}:=\triangle(z_{1},z_{2},w)\cap X$ and $Y_{2}:=\triangle(z_{1},w,z_{3})\cap X$ are Polish spaces having characteristic system $(1,k_{1})$ and $(1,k_{2})$, respectively, such that $0<k_{1},k_{2}<n$ and $k_{1}+k_{2}=n$. As consequence from the induction hypothesis, we have that
\[
\int\limits_{[z_{1},z_{2},w,z_{1}]} f(z) dz=\int\limits_{[z_{1},w,z_{3},z_{1}]}f(z)dz=\textbf{0},
\]
Now, replacing these values in equation \eqref{eq:final}, we hold
\[
\int\limits_{[z_{1},z_{2},z_{3},z_{1}]} f(z) dz=\boldsymbol{0}.
\]
\qed

%%%%%%%%%%%%%%%%%%%%%%%%%%%%%%%%%%%%%%%%%%%%%%%%%%%%%%%%%%%%%%%%%%%%%%%%%%

%%%%%%%%%%%%%%%%%%%%%%%%%%%%%%%%%%
%%%%%%%%%%%%%%%%%%%%%%%%%%


\begin{thebibliography}{99}


\bib{Col1902}{article}{
  title={The april meeting of the American mathematical society},
  author={Cole, F. N.},
  journal = {Bull. Amer. Math. Soc.},
  volume={5},
  number={9},
  pages={423--430},
  year={June 1902},
}



\bib{Ett}{article}{
 author = {H. J. Ettlinger},
 journal = {Annals of Mathematics},
 number = {3},
 pages = {255--270},
 publisher = {Annals of Mathematics},
 title = {Cauchy's Paper of 1814 on Definite Integrals},
 volume = {23},
 year = {1922}
}

\bib{Goursat1884}{article}{
  title={D{\'e}monstration du th{\'e}or{\`e}me de Cauchy},
  author={Goursat, Edouard},
  journal={Acta Mathematica},
  volume={4},
  number={1},
  pages={197--200},
  year={1884},
  publisher={Springer Netherlands},
}

\bib{Goursat1900}{article}{
 author = {Goursat, Edouard},
 journal = {Transactions of the American Mathematical Society},
 number = {1},
 pages = {14--16},
 publisher = {American Mathematical Society},
 title = {Sur La Definition Generale Des Fonctions Analytiques, D'Apres Cauchy},
 volume = {1},
 year = {1900}
}


\bib{kechris}{book}{
  title={Classical descriptive set theory},
  author={Kechris, Alexander},
  volume={156},
  year={2012},
  publisher={Springer Science \& Business Media}
}


\bib{MaSi}{article}{
   author={Mazurkiewicz, S.},
   author={Sierpi\'nski, W.},
   title={Contribution \`a la topologie des ensembles d\`enombrables},
   journal={Fundamenta Mathematicae},
   volume={1},
   date={1920},
   pages={17--27},
}

\bib{More1900}{article}{
  title={A simple proof of the fundamental Cauchy-Goursat theorem},
  author={Moore, Eliakim Hastings},
  journal={Transactions of the American Mathematical Society},
  volume={1},
  number={4},
  pages={499--506},
  year={1900}
}

\bib{Palka}{book}{
  title={An introduction to complex function theory},
  author={Palka, Bruce P.},
  year={1991},
  publisher={Springer Science \& Business Media},
}

\bib{Pring1901}{article}{
  title={Ueber den Goursat’schen Beweis des Cauchy’schen Integralsatzes},
  author={Pringsheim, Alfred},
  journal={Transactions of the American Mathematical Society},
  volume={2},
  number={4},
  pages={413--421},
  year={1901},
}

\bib{Scott1978}{book}{
  title={Cauchy Integral Theorem: a Historical Development of Its Proof},
  author={Scott, Anne Elizabeth},
  series={Master Dissertation},
  year={1978},
  publisher={Oklahoma State University},
}



\end{thebibliography}
\end{document}